\documentclass{article}%
\usepackage{amsfonts}
\usepackage{amsmath}%
\usepackage{amssymb}%
\usepackage{graphicx}
\newtheorem{theorem}{Theorem}

\newtheorem{lemma}[theorem]{Lemma}

\newenvironment{proof}[1][Proof]{\noindent\textbf{#1.} }{\ \rule{0.5em}{0.5em}}
\begin{document}

\title{Intersection Local Time for two Independent Fractional Brownian Motions}
\author{David Nualart
\and Department of Mathematics, University of Kansas
\and 405 Snow Hall, Lawrence, 66045 KS, USA
\and nualart@math.ku.edu, http://www.math.ku.edu/\symbol{126}nualart/
\and Salvador Ortiz-Latorre
\and Facultat de Matem\`{a}tiques, Universitat de Barcelona
\and Gran Via 585 08007 Barcelona, Spain
\and \ \ \ \ \ \ \ \ \ \ \ sortiz@ub.edu \ \ \ \ \ \ \ \ \ \ \ }
\date{}
\maketitle

\begin{abstract}
Let $B^{H}$ and $\widetilde{B}^{H}$ be two independent, $d$-dimensional
fractional Brownian motions with Hurst parameter $H\in\left(  0,1\right)  .$
Assume $d\geq2.$ We prove that the intersection local time of $B^{H}$ and
$\widetilde{B}^{H}$%
\[
I(B^{H},\widetilde{B}^{H})=\int_{0}^{T}\int_{0}^{T}\delta(B_{t}^{H}%
-\widetilde{B}_{s}^{H})dsdt
\]
exists in $L^{2}$ if and only if $Hd<2.$


\textbf{Mathematics Subject Classification MSC2000: }60G15, 60F25, 60G18, 60J55.

\end{abstract}

\section{Introduction}

We consider two independent fractional Brownian motions on $\mathbb{R}%
^{d},d\geq2,$ with the same Hurst parameter $H\in\left(  0,1\right)  .$ This
means that we have two $d$-dimensional independent centered Gaussian processes
$B^{H}=\left\{  B_{t}^{H},t\geq0\right\}  $ and $\widetilde{B}^{H}%
=\{\widetilde{B}_{t}^{H},t\geq0\}$ with covariance structure given by%
\[
\mathbb{E}[B_{t}^{H,i}B_{s}^{H,j}]=\mathbb{E}[\widetilde{B}_{t}^{H,i}%
\widetilde{B}_{s}^{H,j}]=\delta_{ij}R_{H}\left(  s,t\right)  ,
\]
where $i,j=1,...,d,$ $s,t\geq0$ and
\[
R_{H}\left(  s,t\right)  \equiv\frac{1}{2}\left(  t^{2H}+s^{2H}-\left\vert
t-s\right\vert ^{2H}\right)  .
\]
The object of study in this paper will be the intersection local time of
$B^{H}$ and $\widetilde{B}^{H},$ which is formally defined as%
\[
I(B^{H},\widetilde{B}^{H})\equiv\int_{0}^{T}\int_{0}^{T}\delta_{0}(B_{t}%
^{H}-\widetilde{B}_{s}^{H})dsdt,
\]
where $\delta_{0}\left(  x\right)  $ is the Dirac delta function. It is a
measure of the amount of time that the trajectories of the two processes,
$B^{H}$ and $\widetilde{B}^{H},$ intersect on the time interval $\left[
0,T\right]  .$ As we pointed out before, this definition is only formal. In
order to give a rigorous meaning to $I(B^{H},\widetilde{B}^{H})$ we
approximate the Dirac function by the heat kernel
\[
p_{\varepsilon}\left(  x\right)  =\left(  2\pi\varepsilon\right)  ^{-d/2}%
\exp(-\left\vert x\right\vert ^{2}/2\varepsilon),
\]
in $\mathbb{R}^{d}.$ Then, we can consider the following family of random
variables indexed by $\varepsilon>0$
\[
I_{\varepsilon}(B^{H},\widetilde{B}^{H})\equiv\int_{0}^{T}\int_{0}%
^{T}p_{\varepsilon}(B_{t}^{H}-\widetilde{B}_{s}^{H})dsdt,
\]
that we will call the approximated intersection local time of $B^{H}$ and
$\widetilde{B}^{H}$. We are interested in the $L^{2}\left(  \Omega\right)  $
convergence of $I_{\varepsilon}(B^{H},\widetilde{B}^{H})$ as $\varepsilon$
tends to zero.

For $H=1/2,$ the processes $B^{H}$ and $\widetilde{B}^{H}$ are classical
Brownian motions. The intersection local time of independent Brownian motions
has been studied by several authors (see Wolpert \cite{Wo78a} and Geman,
Horowitz and Rosen \cite{GeHoRo84}). The approach of these papers rely on the
fact that the intersection local time of independent Brownian motions can be
seen as the local time at zero of some Gaussian vector field. This approach
easily allows to consider the intersection of $k$ independent Wiener
processes, $k\geq2$. The applications of the intersection local time theory
for Brownian motions range from the construction of relativistic quantum
fields, see Wolpert \cite{Wo78b}, to the construction of the self-intersection
local time for the Brownian motion, see LeGall \cite{LeGall85}. Further
research has been done in order to study such problems for other types of
stochastic processes, mainly L\'{e}vy processes with a particular structure
(strongly symmetric), see Marcus and Rosen \cite{MaRo99}.

In the general case, that is $H\neq1/2,$ only the self-intersection local time
has been studied. Rosen studied in \cite{Ro87} the planar case and a recent
paper by Hu and Nualart \cite{HuNu05} gives a complete picture for the
multidimensional case. On the other hand, Nualart et al. \cite{NuRoTi03} used
a weighted version of the 3-dimensional self-intersection local time for the
study of probabilistic models for vortex filaments based on the fractional
Brownian motion . In recent years the fBm has become an object of intense
study. A stochastic calculus with respect to this process has been developed
by many authors, see Nualart \cite{Nu03} for an extensive account on this
subject. Because of its interesting properties, such as short/long range
dependence and selfsimilarity, the fBm it's being widely used in a variety of
areas such finance, hydrology and telecommunications engineering, see
\cite{NuRoTi03}. Therefore, it seems interesting to study the intersection
local time for this kind of processes.

The aim of this paper is to prove the existence of the intersection local time
of $B^{H}$ and $\widetilde{B}^{H},$ for an $H\neq1/2$ and $d\geq2.$ We have
obtained the following result.

\begin{theorem}
\begin{enumerate}
\item[$(i)$] If $Hd<2,$ then the family of random variables $I_{\varepsilon
}(B^{H},\widetilde{B}^{H})$ converges in $L^{2}\left(  \Omega\right)  $. We
will denote this limit by $I(B^{H},\widetilde{B}^{H}).$

\item[$(ii)$] If $Hd\geq2,$ then
\[
\lim_{\varepsilon\downarrow0}\mathbb{E}[I_{\varepsilon}(B^{H},\widetilde
{B}^{H})]=+\infty
\]
and%
\[
\lim_{\varepsilon\downarrow0}\mathrm{Var}[I_{\varepsilon}(B^{H},\widetilde
{B}^{H})]=+\infty.
\]

\end{enumerate}
\end{theorem}

If $\{B_{t}^{1/2},t\geq0\}$ is a planar Brownian motion, then%
\[
I_{\varepsilon}=\int_{0}^{T}\int_{0}^{T}\delta_{0}\left(  B_{s}^{1/2}%
-B_{t}^{1/2}\right)  dsdt
\]
diverges almost sure, when $\varepsilon$ tends to zero. Varadhan, in
\cite{Va69}, proved that the renormalized self-intersection local time defined
as $\lim_{\varepsilon\rightarrow0}(I_{\varepsilon}-\mathbb{E}[I_{\varepsilon
}]),$ exists in $L^{2}\left(  \Omega\right)  $. Condition $\left(  ii\right)
$ implies that Varadhan renormalization does not converge in this case.

For $Hd\geq2,$ according to the previous theorem, $I_{\varepsilon}%
(B^{H},\widetilde{B}^{H})$ doesn't converge in $L^{2}\left(  \Omega\right)  $
and therefore $I(B^{H},\widetilde{B}^{H}),$ the intersection local time of
$B^{H} $ and $\widetilde{B}^{H},$ doesn't exist. The proof of Theorem 1.1 rest
on Lemma \ref{LemaDetFinite}, which deals with the integral of a negative
power of the determinant of some covariance matrix.

The paper is organized as follows. In Section 2 we prove Theorem 1.1. In order
to clarify the exposition, some technical lemmas needed in the proof are
stated and proved in the Appendix.

\section{Intersection Local Time of $B^{H}$ and $\widetilde{B}^{H},$ Case
$Hd<2$}

Let $B^{H}$ and $\widetilde{B}^{H}$ two independent fractional Brownian
motions on $\mathbb{R}^{d}$ with the same Hurst parameter $H\in\left(
0,1\right)  .$

Using the following classical equality
\[
p_{\varepsilon}\left(  x\right)  =\frac{1}{\left(  2\pi\right)  ^{d}}%
\int_{\mathbb{R}^{d}}e^{i\left\langle \xi,x\right\rangle }e^{-\varepsilon
\frac{\left\vert \xi\right\vert ^{2}}{2}}d\xi
\]
from Fourier analysis, and the definition of $I_{\varepsilon}(B^{H}%
,\widetilde{B}^{H}),$ we obtain
\begin{equation}
I_{\varepsilon}(B^{H},\widetilde{B}^{H})=\frac{1}{\left(  2\pi\right)  ^{d}%
}\int_{0}^{T}\int_{0}^{T}\int_{\mathbb{R}^{d}}e^{i\langle\xi,B_{t}%
^{H}-\widetilde{B}_{s}^{H}\rangle}e^{-\varepsilon\frac{\left\vert
\xi\right\vert ^{2}}{2}}d\xi dsdt.\label{FourierIEps}%
\end{equation}
Therefore,%

\begin{align}
\mathbb{E}[I_{\varepsilon}(B^{H},\widetilde{B}^{H})]  & =\frac{1}{\left(
2\pi\right)  ^{d}}\int_{0}^{T}\int_{0}^{T}\int_{\mathbb{R}^{d}}\mathbb{E[}%
e^{i\langle\xi,B_{t}^{H}-\widetilde{B}_{s}^{H}\rangle}]e^{-\varepsilon
\frac{\left\vert \xi\right\vert ^{2}}{2}}d\xi dsdt\nonumber\\
& =\frac{1}{\left(  2\pi\right)  ^{d}}\int_{0}^{T}\int_{0}^{T}\int
_{\mathbb{R}^{d}}e^{-(\varepsilon+s^{2H}+t^{2H})\frac{\left\vert
\xi\right\vert ^{2}}{2}}d\xi dsdt\nonumber\\
& =\frac{1}{\left(  2\pi\right)  ^{d/2}}\int_{0}^{T}\int_{0}^{T}%
(\varepsilon+s^{2H}+t^{2H})^{-d/2}dsdt,\label{ExpIEps}%
\end{align}
where we have used that $\langle\xi,B_{t}^{H}-\widetilde{B}_{s}^{H}\rangle\sim
N(0,\left\vert \xi\right\vert ^{2}\left(  s^{2H}+t^{2H}\right)  ),$ so
\[
\mathbb{E[}e^{i\langle\xi,B_{t}^{H}-\widetilde{B}_{s}^{H}\rangle}%
]=e^{-(s^{2H}+t^{2H})\frac{\left\vert \xi\right\vert ^{2}}{2}},
\]
and the fact that
\[
\int_{\mathbb{R}^{d}}e^{-(\varepsilon+s^{2H}+t^{2H})\frac{\left\vert
\xi\right\vert ^{2}}{2}}d\xi=\left(  \frac{2\pi}{\varepsilon+s^{2H}+t^{2H}%
}\right)  ^{d/2}.
\]

According to the representation $\left(  \ref{FourierIEps}\right)  $ for
$I_{\varepsilon}(B^{H},\widetilde{B}^{H}),$ we have that
\begin{align}
\mathbb{E}[I_{\varepsilon}^{2}(B^{H},\widetilde{B}^{H})]  & =\frac{1}{\left(
2\pi\right)  ^{2d}}\int_{\left[  0,T\right]  ^{4}}\int_{\mathbb{R}^{2d}%
}\mathbb{E}[e^{i(\langle\xi,B_{t}^{H}-\widetilde{B}_{s}^{H}\rangle+\langle
\eta,B_{v}^{H}-\widetilde{B}_{u}^{H}\rangle)}]\nonumber\\
& \times e^{-\varepsilon\frac{\left\vert \xi\right\vert ^{2}+\left\vert
\eta\right\vert ^{2}}{2}}d\xi d\eta dsdtdudv.\label{ExpIEpsSqr}%
\end{align}
Let introduce some notation that we will use throughout this paper,%
\begin{align*}
\lambda & =\lambda\left(  s,t\right)  =s^{2H}+t^{2H},\\
\rho & =\rho\left(  u,v\right)  =u^{2H}+v^{2H},
\end{align*}
and
\[
\mu=\mu\left(  s,t,u,v\right)  =\frac{1}{2}\left(  s^{2H}+t^{2H}+u^{2H}%
+v^{2H}-\left\vert t-v\right\vert ^{2H}-\left\vert s-u\right\vert
^{2H}\right)  .
\]
Notice that $\lambda$ is the variance of $B_{t}^{H,1}-B_{s}^{H,2},$ $\rho$ is
the variance of $B_{v}^{H,1}-B_{u}^{H,2}$ and $\mu$ is the covariance between
$B_{t}^{H,1}-B_{s}^{H,2}$ and $B_{v}^{H,1}-B_{u}^{H,2}$, where $B^{H,1}$ and
$B^{H,2}$ are independent one-dimensional fractional Brownian motions with
Hurst parameter $H.$

Using that $\langle\xi,B_{t}^{H}-\widetilde{B}_{s}^{H}\rangle+\langle
\eta,B_{v}^{H}-\widetilde{B}_{u}^{H}\rangle\sim N(0,\lambda\left\vert
\xi\right\vert ^{2}+\rho\left\vert \eta\right\vert ^{2}+2\mu\langle\xi
,\eta\rangle)$ and $\left(  \ref{ExpIEpsSqr}\right)  $ we can write for all
$\varepsilon>0$%
\begin{align}
& \mathbb{E}[I_{\varepsilon}^{2}(B^{H},\widetilde{B}^{H})]\nonumber\\
& =\frac{1}{\left(  2\pi\right)  ^{2d}}\int_{\left[  0,T\right]  ^{4}}%
\int_{\mathbb{R}^{2d}}e^{-\frac{1}{2}\{(\lambda+\varepsilon)\left\vert
\xi\right\vert ^{2}+(\rho+\varepsilon)\left\vert \eta\right\vert ^{2}%
+2\mu\langle\xi,\eta\rangle\}}d\xi d\eta dsdtdudv\nonumber\\
& =\frac{1}{\left(  2\pi\right)  ^{d}}\int_{\left[  0,T\right]  ^{4}}%
((\lambda+\varepsilon)(\rho+\varepsilon)-\mu^{2})^{-d/2}%
dsdtdudv.\label{ExpIEpsSqrFinal}%
\end{align}
The last equality follows from the well known fact that
\[
\int_{\mathbb{R}^{2d}}e^{-\frac{1}{2}\langle x,Ax\rangle}dx=\frac{\left(
2\pi\right)  ^{d}}{\left(  \det A\right)  ^{1/2}},
\]
with%
\[
A=Id_{d}\otimes%
\begin{pmatrix}
\lambda+\varepsilon & \mu\\
\mu & \rho+\varepsilon
\end{pmatrix}
,
\]
where $Id_{d}$ is the $d$-dimensional identity matrix and $\otimes$ denotes
the Kronecker product of matrices. We also have that
\begin{align*}
\det A  & =\det\left(  Id_{d}\otimes%
\begin{pmatrix}
\lambda+\varepsilon & \mu\\
\mu & \rho+\varepsilon
\end{pmatrix}
\right)  \\
& =\left(  \det\left(  Id_{d}\right)  \right)  ^{2}\cdot\left(  \det%
\begin{pmatrix}
\lambda+\varepsilon & \mu\\
\mu & \rho+\varepsilon
\end{pmatrix}
\right)  ^{d}\\
& =((\lambda+\varepsilon)(\rho+\varepsilon)-\mu^{2})^{d}.
\end{align*}

\begin{proof}
[\textbf{Proof of Theorem 1}]Suppose first that $Hd<2.$ A slight extension of
$\left(  \ref{ExpIEpsSqrFinal}\right)  $ yields%
\[
\mathbb{E}[I_{\varepsilon}(B^{H},\widetilde{B}^{H})I_{\eta}(B^{H}%
,\widetilde{B}^{H})]=\int_{\left[  0,T\right]  ^{4}}((\lambda+\varepsilon
)(\rho+\eta)-\mu^{2})^{-d/2}dsdtdudv.
\]
Consequently, a necessary and sufficient condition for the convergence in
$L^{2}\left(  \Omega\right)  $ of $I_{\varepsilon}(B^{H},\widetilde{B}^{H})$
is that
\[
\int_{\left[  0,T\right]  ^{4}}(\lambda\rho-\mu^{2})^{-d/2}dsdtdudv<+\infty.
\]
Then the result follows from Lemma \ref{LemaDetFinite}.

Now suppose that $Hd\geq2,$ then from $\left(  \ref{ExpIEps}\right)  $ and
using monotone convergence theorem%
\[
\lim_{\varepsilon\downarrow0}\mathbb{E}[I_{\varepsilon}(B^{H},\widetilde
{B}^{H})]=\int_{0}^{T}\int_{0}^{T}\left(  s^{2H}+t^{2H}\right)  ^{-d/2}dsdt,
\]
and this integral is divergent by Lemma \ref{LemaIntegralLambda}. According to
the expression $\left(  \ref{ExpIEps}\right)  $ for $\mathbb{E}[I_{\varepsilon
}(B^{H},\widetilde{B}^{H})]$ and the expression $\left(  \ref{ExpIEpsSqrFinal}%
\right)  $ for $\mathbb{E}[I_{\varepsilon}^{2}(B^{H},\widetilde{B}^{H})]$ we
obtain%
\begin{align*}
\lim_{\varepsilon\downarrow0} &  \mathrm{Var}[I_{\varepsilon}(B^{H}%
,\widetilde{B}^{H})]=\lim_{\varepsilon\downarrow0}\left\{  \mathbb{E}%
[I_{\varepsilon}(B^{H},\widetilde{B}^{H})^{2}]-\left(  \mathbb{E}%
[I_{\varepsilon}(B^{H},\widetilde{B}^{H})]\right)  ^{2}\right\} \\
&  =\int_{\left[  0,T\right]  ^{4}}(\lambda\rho-\mu^{2})^{-d/2}-\left(
\lambda\rho\right)  ^{-d/2}dsdtdudv.
\end{align*}
Set
\begin{equation}
D_{\varepsilon}:=\{(s,t,u,v)\in\mathbb{R}_{+}^{4}\ |\ s^{2}+t^{2}+u^{2}%
+v^{2}\leq\varepsilon^{2}\}.\label{ExpDEps}%
\end{equation}
We can find $\varepsilon>0$ such that $D_{\varepsilon}\subset\left[
0,T\right]  ^{4}.$ Making a change to spherical coordinates, as the integrand
is always positive, we have
\begin{align*}
& \int_{\left[  0,T\right]  ^{4}}(\lambda\rho-\mu^{2})^{-d/2}-\left(
\lambda\rho\right)  ^{-d/2}dsdtdudv\\
& \geq\int_{D_{\varepsilon}}(\lambda\rho-\mu^{2})^{-d/2}-\left(  \lambda
\rho\right)  ^{-d/2}dsdtdudv=\int_{0}^{\varepsilon}r^{3-2Hd}dr\int_{\Theta
}\Psi\left(  \theta\right)  d\theta,
\end{align*}
where the integral in $r$ is convergent if and only if $Hd<2,$ and the angular
integral is different from zero thanks to the positivity of the integrand.
Therefore, if $Hd\geq2,$ then
\[
\lim_{\varepsilon\downarrow0}\mathrm{Var}[I_{\varepsilon}(B^{H},\widetilde
{B}^{H})]=+\infty.
\]

\end{proof}

\section{Appendix}

For clarity of exposition, we state and prove some technical lemmas in this appendix.

\begin{lemma}
\label{LemaPsi}Let $\alpha>0,$ and let%
\begin{equation}
\gamma\left(  \alpha,x\right)  \equiv\int_{0}^{x}e^{-y}y^{\alpha
-1}dy\label{DefGamma}%
\end{equation}
be the lower incomplete gamma function. Then for all $\varepsilon<\alpha$ and
$x>0,$
\[
\gamma\left(  \alpha,x\right)  \leq K\left(  \alpha\right)  x^{\varepsilon},
\]
where $K\left(  \alpha\right)  \equiv\frac{1}{\alpha}\vee\Gamma\left(
\alpha\right)  $ and $\Gamma\left(  \alpha\right)  =\gamma\left(
\alpha,+\infty\right)  $.
\end{lemma}

\begin{proof}
If $x\geq1,$
\[
\gamma\left(  \alpha,x\right)  \leq\Gamma\left(  \alpha\right)  x^{\varepsilon
},
\]
for all $\varepsilon>0.$ On the other hand, if $x<1,$%
\[
\gamma\left(  \alpha,x\right)  \leq\int_{0}^{x}y^{\alpha-1}dy=\frac{x^{\alpha
}}{\alpha}=\frac{1}{\alpha}x^{\varepsilon},
\]
if $\varepsilon<\alpha.$
\end{proof}

\begin{lemma}
\label{LemaIntegralLambda}The following integral
\[
\int_{0}^{T}\int_{0}^{T}\left(  s^{2H}+t^{2H}\right)  ^{-d/2}dsdt,
\]
is finite if and only if $Hd<2.$
\end{lemma}

\begin{proof}
It easily follows from a polar change of coordinates.
\end{proof}

\begin{lemma}
\label{LemaDetFinite}Let
\[
A_{T}\equiv\int_{\left[  0,T\right]  ^{4}}(\lambda\rho-\mu^{2})^{-d/2}%
dsdtdudv,
\]
then $A_{T}<+\infty$ if and only if $Hd<2.$
\end{lemma}

\begin{proof}
The necessary condition follows from a spherical change of coordinates. We can
find $\varepsilon>0$ such that $D_{\varepsilon}\subset\left[  0,T\right]
^{4},$ where $D_{\varepsilon}$ is given in $\left(  \ref{ExpDEps}\right)  .$
As the integrand in $A_{T}$ is always positive we have
\[
A_{T}\geq\int_{D_{\varepsilon}}(\lambda\rho-\mu^{2})^{-d/2}dsdtdudv=\int
_{0}^{\varepsilon}r^{3-2Hd}dr\int_{\Theta}\phi\left(  \theta\right)  d\theta,
\]
where the integral in $r$ is convergent if and only if $Hd<2,$ and the angular
integral is different from zero thanks to the positivity of the integrand.
Therefore, if $Hd\geq2,$ then $A_{T}=+\infty.$

Suppose now that $Hd<2.$ By symmetry we have that
\[
A_{T}=4\int_{\mathcal{T}}\left(  \lambda\rho-\mu^{2}\right)  ^{-d/2}dsdtdudv,
\]
where
\[
\mathcal{T}\equiv\{\left(  s,t,u,v):0<v<t,0<t\leq T,0<u<s,0<s\leq T\right)
\}.
\]
Notice that
\[
\lambda\rho-\mu^{2}=\det\mathrm{Var}\left(  Z\right)  ,
\]
where $Z\equiv(B_{t}^{H,1}-\widetilde{B}_{s}^{H,1},B_{v}^{H,1}-\widetilde
{B}_{u}^{H,1}).$ Due to the independence of $B^{H}$ and $\widetilde{B}^{H},$
we have that
\[
\mathrm{Var}\left(  Z\right)  =\mathrm{Var}(B_{t}^{H,1},B_{v}^{H,1}%
)+\mathrm{Var}(\widetilde{B}_{s}^{H,1},\widetilde{B}_{u}^{H,1}),
\]
and
\[
\lambda\rho-\mu^{2}\geq\det(\mathrm{Var}(B_{t}^{H,1},B_{v}^{H,1}%
))+\det(\mathrm{Var}(\widetilde{B}_{s}^{H,1},\widetilde{B}_{u}^{H,1})),
\]
because the matrices $\mathrm{Var}(B_{t}^{H,1},B_{v}^{H,1})$ and
$\mathrm{Var}(\widetilde{B}_{s}^{H,1},\widetilde{B}_{u}^{H,1})$ are strictly
positive definite (see A8, $\left(  \mathrm{viii}\right)  $ in \cite{Mu82}).
Then
\[
A_{T}\leq4\int_{\mathcal{T}}(\varphi\left(  t,v\right)  +\varphi\left(
s,u\right)  )^{-d/2}dsdtdudv,
\]
where
\[
\varphi\left(  t,v\right)  \equiv\det(\mathrm{Var}(B_{t}^{H,1},B_{v}%
^{H,1}))=t^{2H}v^{2H}-\frac{1}{4}\left(  t^{2H}+v^{2H}-\left\vert
t-v\right\vert ^{2H}\right)  ^{2}.
\]
Using Fubini's Theorem and
\[
\lambda^{-\alpha}=\frac{1}{\Gamma\left(  \alpha\right)  }\int_{0}^{+\infty
}e^{-\lambda z}z^{\alpha-1}dz,
\]
for all $\lambda,\alpha>0$, we obtain
\begin{align}
& \int_{\mathcal{T}}(\varphi\left(  t,v\right)  +\varphi\left(  s,u\right)
)^{-d/2}dsdtdudv\nonumber\\
& =\frac{1}{\Gamma(\frac{d}{2})}\int_{\mathcal{T}}\int_{0}^{+\infty
}e^{-(\varphi\left(  t,v\right)  +\varphi\left(  s,u\right)  )z}z^{\frac{d}%
{2}-1}dzdsdtdudv\nonumber\\
& =\frac{1}{\Gamma(\frac{d}{2})}\int_{0}^{+\infty}z^{\frac{d}{2}-1}%
A^{2}\left(  z\right)  dz,\label{IntegralA}%
\end{align}
where
\[
A\left(  z\right)  \equiv\int_{0}^{T}\int_{0}^{t}e^{-\varphi\left(
t,v\right)  z}dvdt.
\]
As $A\left(  z\right)  <+\infty,$ for all $z\in\lbrack0,1],$ the integral
$\left(  \ref{IntegralA}\right)  $ is convergent in a neighborhood of zero.
Hence, we have to study the convergence of
\[
\int_{1}^{+\infty}z^{\frac{d}{2}-1}A^{2}\left(  z\right)  dz.
\]
Due to the homogeneity of order $4H$ of $\varphi\left(  t,v\right)  ,$ if we
make the change of coordinates $t=z^{-\frac{1}{4H}}x,v=z^{-\frac{1}{4H}}y,$ we
obtain%
\[
\int_{1}^{+\infty}z^{\frac{d}{2}-1}A^{2}\left(  z\right)  dz=\int_{1}%
^{+\infty}z^{\frac{d}{2}-1-\frac{1}{H}}\left(  \int_{0}^{Tz^{\frac{1}{4H}}%
}\int_{0}^{x}e^{-\varphi\left(  x,y\right)  }dydx\right)  ^{2}dz.
\]
Now, using that $\{(x,y):0<x<Tz^{\frac{1}{4H}},0<y<x\}\subset\{(x,y):x^{2}%
+y^{2}\leq2T^{2}z^{\frac{1}{2H}}\}\equiv S,$ and making a polar change of
coordinates we have%
\[
\int_{0}^{Tz^{\frac{1}{4H}}}\int_{0}^{x}e^{-\varphi\left(  x,y\right)
}dydx\leq\int_{S}e^{-\varphi\left(  x,y\right)  }dydx=\int_{0}^{\pi/4}\int
_{0}^{\sqrt{2}Tz^{\frac{1}{4H}}}re^{-r^{4H}\varphi\left(  \theta\right)
}drd\theta,
\]
where $\varphi\left(  \theta\right)  \equiv\varphi\left(  \cos\theta
,\sin\theta\right)  .$ After the new change of variable $x=r^{4H}%
\varphi\left(  \theta\right)  ,$ the last integral is equal to
\begin{align*}
& \int_{0}^{\pi/4}\varphi\left(  \theta\right)  ^{-\frac{1}{2H}}\int
_{0}^{2^{2H}T^{4H}z\varphi\left(  \theta\right)  }\frac{1}{4H}e^{-x}%
x^{\frac{1}{2H}-1}dxd\theta\\
& =\frac{1}{4H}\int_{0}^{\pi/4}\varphi\left(  \theta\right)  ^{-\frac{1}{2H}%
}\gamma\left(  (2H)^{-1},2^{2H}T^{4H}z\varphi\left(  \theta\right)  \right)
d\theta,
\end{align*}
where $\gamma\left(  \alpha,x\right)  $ is given by $\left(  \ref{DefGamma}%
\right)  .$ Applying Lemma \ref{LemaPsi},%
\begin{align*}
& \int_{1}^{+\infty}z^{\frac{d}{2}-1}A^{2}\left(  z\right)  dz\\
& \leq\frac{1}{16H^{2}}\int_{1}^{+\infty}z^{\frac{d}{2}-1-\frac{1}{H}}\left(
\int_{0}^{\pi/4}\varphi\left(  \theta\right)  ^{-\frac{1}{2H}}\gamma\left(
(2H)^{-1},2^{2H}T^{4H}z\varphi\left(  \theta\right)  \right)  d\theta\right)
^{2}dz\\
& \leq\frac{2^{4H\varepsilon}T^{8H\varepsilon}}{16H^{2}}K^{2}\left(  \frac
{1}{2H}\right)  \int_{1}^{+\infty}z^{\frac{d}{2}-1-\frac{1}{H}+2\varepsilon
}dz\left(  \int_{0}^{\pi/4}\varphi\left(  \theta\right)  ^{\varepsilon
-\frac{1}{2H}}d\theta\right)  ^{2}.
\end{align*}
The integral in $z$ is convergent provided $\varepsilon<\frac{2-Hd}{4H}.$ It's
an exercise of computation of limits to prove that $\varphi\left(
\theta\right)  \sim\theta^{2H}$ as $\theta\downarrow0$ and $\varphi\left(
\theta\right)  \sim(\pi/4-\theta)^{2H}$ as $\theta\uparrow\pi/4, $ the main
tool is to substitute the trigonometric functions by their first order
approximations at the respective points$.$ As a consequence, the integral
\[
\int_{0}^{\pi/4}\varphi\left(  \theta\right)  ^{\varepsilon-\frac{1}{2H}%
}d\theta
\]
is always convergent.
\end{proof}

\end{document}